\documentclass[10pt]{amsart}

\usepackage{amsmath}
\usepackage{amssymb}
\usepackage{url}
\usepackage{verbatim}
\DeclareMathOperator{\lcm}{lcm}
\newcommand{\nc}{\newcommand}
\nc{\gp}{\texttt{gp}}
\nc{\Q}{{\bf Q}}
\nc{\Z}{{\bf Z}}
\nc{\Qp}{\Q_p}
\nc{\Zp}{\Z_p}
\nc{\C}{{\bf C}}
\nc{\bbP}{{\bf P}}
\nc{\R}{{\bf R}}
\nc{\F}{{\bf F}}

\nc{\Qbar}{\overline{\Q}}
\nc{\un}[1]{#1^{\mathrm{un}}}
\nc{\emid}{\mid\kern -0.1ex\mid}

\nc{\OK}{{\mathcal{O}}_K}
\nc{\cK}{{\mathcal{K}}}
\nc{\tr}{\mathop{\rm Tr}}
\nc{\disc}{\mathop{\rm disc}}
\nc{\ord}{\mbox{ord}}
\nc{\cond}{\mathop{\rm cond}}
\nc{\Gal}{\mathop{\rm Gal}}
\nc{\Sym}{\mathop{\rm Sym}}
\nc{\Aut}{\mathop{\rm Aut}}
\nc{\Quot}{\mathop{\rm Quot}}
\nc{\cmmt}[1]{}

\newtheorem{Proposition}[table]{Proposition}
\makeatletter
\@addtoreset{table}{section}
\makeatother

 \title{
 A database of local fields}
 \author{John W.\ Jones and David P.\ Roberts}

\begin{document}
%\begin{frontmatter}
%  \title{A Database of Local Fields}
%  \author{John W.\ Jones}
%  \address{Department of Mathematics and Statistics, Arizona State
%    University, Box 871804, Tempe, AZ 85287}
%  \ead{jj@asu.edu}
%  \author{David P.\ Roberts}
%  \address{Division of Science and Mathematics, University of
%      Minnesota--Morris, Morris, MN 56267}
%  \ead{roberts@mrs.umn.edu}
  
\begin{abstract}
  We describe our online database of finite extensions of $\Qp$, and how 
  it can be used to facilitate local analysis of number fields.
\end{abstract}
%\end{frontmatter}
\maketitle
\section{Introduction}

\subsection{Overview}
Given a number field $K$, one has for each prime $p$ its 
associated $p$-adic algebra,
\[K \otimes \Qp \cong \prod_{i=1}^g K_{p,i}\, .\]
Here the $K_{p,i}$ are fields, each a finite extension of $\Qp$.
For investigating some problems about number fields, 
it suffices to know just basic invariants of the
$K_{p,i}$, such as ramification index and residual degree.  For
other investigations, it is essential to have much more refined 
information, such as local Galois groups and slopes measuring wildness
of ramification.  

To facilitate refined analysis of number fields,
we have constructed 
%an 
a
%online
database of $p$-adic fields, available
at 
\url{http://math.asu.edu/~jj/localfields}.
%\url{math.asu.edu/~jj/localfields}.
Let $\cK(p,n)$ be the set
of isomorphism classes of 
degree $n$ extensions of $\Q_{p}$.  The sets $\cK(p,n)$ 
are finite, with general mass formulas counting these fields with 
certain weights being known \cite{serre-mass,krasner,pauli-roblot}.  Our database 
presents some of the sets $\cK(p,n)$ in a complete and
easy-to-use way.   The philosophy behind the database is that
the intricate local considerations needed to construct it should
be done once and then recorded.  Thereafter, a local result
can be obtained by mechanical appeal to the database whenever
it is needed in a global situation.

\subsection{Fields in the database} 
When $n$ is not divisible by $p$, all fields in $\cK(p,n)$
are tame, and so $\cK(p,n)$ is relatively easy to describe.
Our database treats these fields 
%and all other tame fields
dynamically, without restriction on  $p$ or $n$.  
The first case involving wild fields is $n=p$.  
This case is also relatively easy to describe in a way uniform in $p$;
for example, $|\cK(p,p)| = p^{2}+1$ for $p$ odd.  Again, our database
treats these fields without restriction on $p$.

\begin{table}[htb] \caption{\label{cards}The number $|\cK(p,n)|$ of
    isomorphism 
classes of $p$-adic fields of degree $n$, 
for $p<30$ and $n<10$.  The entries corresponding to the
five cases which we treat individually are underlined.} 
\[
\begin{array}{r|rrrrrrrrrr}
n & 2 & 3 & \hspace{10pt} 5 & 
\hspace{10pt} 7 & \hspace{5pt} 11 & \hspace{5pt} 13 & 
\hspace{5pt} 17 & \hspace{5pt} 19 & \hspace{5pt} 23 & \hspace{5pt} 29 \\
\hline
1 & 1 & 1 & 1 & 1 & 1 & 1 & 1 & 1 & 1 & 1 \\
2 & 7 & 3 & 3 & 3 & 3 & 3 & 3 & 3 & 3 & 3 \\
3 & 2 & 10 & 2 & 4 & 2 & 4 & 2 & 4 & 2 & 2  \\
4 &\underline{59} & 5 & 7 & 5& 5 & 7& 7 & 5& 5& 7 \\
5 & 2 & 2 & 26 & 2 & 6 & 2 & 2 & 2 & 2 & 2  \\
6 &\underline{47}&\underline{75} & 7 & 12 & 7& 12& 7& 12& 7& 7 \\
7 & 2 & 2 & 2 & 50 & 2 & 2 & 2 & 2 & 2 & 8 \\
8 & \underline{1823}&8& 11& 8& 8& 11& 15& 8& 8& 11 \\
9 & 3 &\underline{795}& 3& 7& 3& 7& 3& 13& 3& 3 \\
\end{array}
\]
\end{table}

%The most visible part of our database are
The most visible parts of our database are
tables explicitly describing $\cK(p,n)$ for small $p$ and $n$.
%$p=2$, $3$, $5$, $7$ and $n \leq 9$.
The numbers 
$|\cK(p,n)|$ for $p<30$ and $n<10$ 
are listed in Table~\ref{cards}. The table for $\cK(p,n)$ in the database has one line
for each isomorphism class of $p$-adic field of degree $n$ and
gives a defining polynomial for the field and many invariants 
of the field.   
Our tables provide many illustrations of the 
relatively easy cases discussed in the previous paragraph.  However
their main function is to cover the five harder cases with $n < 10$, 
namely $(p,n)=(2,4)$, 
$(2,6)$, $(2,8)$, $(3,6)$, and $(3,9)$.

\subsection{Sections of this paper}
Section~\ref{construction} discusses how we found our 
lists of defining polynomials.  It treats first the 
tame and $n=p$ cases systematically, and then describes our {\em ad 
hoc} approach to the five harder cases.  All parts of
the $2$-adic quartic cases have received detailed attention previously, 
for example \cite{weil} for the one $A_{4}$ and the three $S_{4}$ 
extensions and \cite{naito} for the thirty-six $D_{4}$ extensions.  
However even the mere listing of defining polynomials constitutes
a new result in the remaining
four cases.   For example, while some classes of 
octic $2$-adic fields have been studied completely
by others \cite{weil,bayer-rio}, these fields represent only 
a small subset of the full set of 1823 octic 2-adic fields.  

    Section~\ref{invariants}
discusses how we computed the invariants for each field.
We restrict discussion of details to the  cases $n \leq 7$.  
The analogous details for $2$-adic octics and $3$-adic 
nonics are very much more intricate. We are treating these cases
in separate papers, each of which will 
include sample applications to number fields
%, 
of the same degree, where the same
Galois theory applies.

%Section~? discusses how we found our complete lists of 
%defining polynomials.  It treats first the relatively easy
%unramified, tamely ramified, and $p=n=e$ cases, and
%then the general case.  

%   Section~\ref{sampleuse} illustrates
%how our database can be used in practice.  
%It describes the two interactive features of our
%database by means of a global example involving degree 13 
%polynomials.  

%and explains the two interactive features of our database.  
%These each take a monic polynomial $f(x) \in \Z[x]$ 
%as input.  Importantly, the programs work with the local fields
%on our tables and also also tame and degree $p$ fields.  
%So they give complete results in many global situations where
%the degree of the input $f(x)$ is well beyond $9$. 
 
   Section~4 describes the two interactive features of our
database, what we call the $p$-adic identifier and the
Galois root discriminant calculator.  These are designed to 
maximize the utility of our database for applications.  
One application we have in mind is to assist in matching number fields 
to automorphic forms of various sorts, as here complete understanding
of ramification is very useful.  Section~4 presents another 
application, one that stays within the confines of traditional
algebraic number theory.

\section{A complete irredundant list of defining polynomials}
\label{construction}
In this section, we describe how we chose the polynomials defining the 
fields in the database. 
Sections \ref{unram}-\ref{sec:degree-p-ramified} deal
with unramified, tamely ramified, and degree $p$ extensions of 
$\Q_p$, respectively.  
%In those cases, the extensions are well understood so a computer can
%generate the data it needs on the fly.  For the sake of consistency,
%choices have to be made of a particular defining polynomial for each
%extension.
%Finally, 
Section~\ref{sec:wild-extens-comp} deals with the remaining cases
--- wildly ramified extensions where the degree is composite.

\subsection{Unramified extensions}
\label{unram}
Unramified extensions of $\Qp$ are very simple, there being a unique
one for each degree $n$, up to isomorphism.   The only
task is to choose a defining polynomial for each.  In the sequel, we 
will sometimes drop qualifiers like ``up to isomorphism,'' as they are
always present and our meaning is clear.  

Since the extension of degree $n$ of $\Qp$ corresponds to the unique
degree $n$ extension of the residue field, one option is to use
Conway polynomials for these extensions since they are
considered the standard choices for defining $\F_{p^n}$ over $\F_p$.
However, Conway
polynomials can be expensive to compute, primarily because they
are required to satisfy a compatibility condition which is not used
here.
% Since unramified extensions
% are so thoroughly understood, it does not make sense to pick defining
% polynomials for them which are difficult to compute.

Instead, we pick defining polynomials which are in the same spirit,
but with fewer restrictions.  We compute the ``first'' polynomial over
$\F_p$ which has roots which are primitive, i.e., of multiplicative
order $p^n-1$.
Here we use the same lexicographic ordering as for Conway
polynomials.  We write polynomials in the form
$f(x) = x^n - a_{n-1} x^{n-1} + a_{n-2} x^{n-2}
-\cdots$ and $g(x) = x^n - b_{n-1} x^{n-1} + b_{n-2} x^{n-2}-\cdots$
with $a_i$ and $b_i$ between $0$ and $p-1$ inclusive.  Then we
define $f < g$ iff there exists $k$ with $a_i = b_i$ for all $i> k$
and $a_k<b_k$.  This normalization also defines how we will represent
the polynomials in $\Z[x]\subset\Qp[x]$.
Note that for defining $\Qp$ itself, our choice leads to the ``degree
1 Conway polynomial'' $x-r$, where $r$ is the first primitive root
modulo $p$.

\subsection{Tame extensions}

Our starting point is the following standard result on totally
ramified tame extensions.  The statement here is a special case of 
Theorem~7.2 in \cite{pauli-roblot}.

\begin{Proposition}\label{tame-thm}
  Let $K^{u}$ be an unramified extension of $\Q_p$ with degree $f$.
Let $\zeta\in K^{u}$ be a primitive $(p^f-1)^{\rm st}$ root of unity.
Let $e$ be a positive integer with $p \nmid
e$.
\begin{enumerate}
\item The totally ramified degree $e$ extensions of $K^{u}$, are given
by roots of polynomials $h_{e,r}(x) = x^e - \zeta^r p$.
\item Two
such polynomials $h_{e,r}$ and $h_{e,r'}$ yield $K^{u}$-isomorphic extensions
iff $r\equiv r' \pmod{ \gcd(e,p^f-1 )}$.
\item If a monic polynomial
$g$ satisfies $g\equiv h_{e,r} \pmod{p^2}$, then $g$ defines the same
extension as $h_{e,r}$.
\end{enumerate}
\end{Proposition}

To apply the proposition, we take $K^{u} = \Q_{p}[\alpha]/h(\alpha)$ 
where $h$ is the degree $f$ polynomial chosen in the previous 
subsection.  We consider $x^{e} - \alpha^{r} p$, as the third part of
Proposition~\ref{tame-thm} lets us replace $\zeta$ by $\alpha$.

Since $\Gal(K^{u}/\Qp)$ is generated by Frobenius $\sigma$, with
\[ \sigma(\alpha) \equiv \alpha^p \pmod{p},\]
the polynomials $x^e -\alpha^r p$ give conjugate
extensions for $r$ which differ multiplicatively by a power of $p$.
Taking the norm of $x^e-\alpha^r p$ to $\Qp$, we get an irreducible
polynomial iff the orbit of $r$ in $\Z/(p^f-1)\Z$ under multiplication
by $p$ has length $f$.

%It is not hard to show that under these
%conditions, the resulting polynomial is also irreducible over $\Qp$.
%and is independent of the choice of the root $\alpha$.

%isomorphic extensions over $\Qp(\alpha)$ for each
%congruence class of $r$ modulo $\gcd(e, p^f-1)$, and give conjugate
%extensions for $r$ which differ multiplicatively by a power of $p$.

Our recipe for picking defining polynomials of tamely ramified
extensions with given $e$ and $f$ is as follows.  Let
$g=\gcd(e,p^f-1)$ and partition $\Z/g\Z$ into orbits under
multiplication by $p$.  These will correspond to the desired
extensions of $\Qp$.  For each orbit $\mathcal{O}\subseteq \Z/g\Z$, 
we lift its
elements to $\Z/(p^f-1)\Z$ and consider them under multiplication by
$p$.  If there is an orbit of length $f$, take
the smallest $r\ge 0$ contained in such an orbit.  Then the norm of
$x^e-\alpha^r p$ to $\Qp$ will be irreducible.
If there are no lifts to an orbit of length $f$ for our orbit
$\mathcal{O}$, we apply the following ``root shift'' procedure. 
We take the smallest $r \ge 0$ representing an element
of the orbit and consider polynomials $(x+k\alpha)^e - \alpha^r p$
with $k = 1$, $2$, $3$, \dots  and take their norms to $\Qp[x]$.
The first norm which is irreducible is our preferred
defining polynomial.

For example, to generate the sextic tame extensions of $\Q_5$ with
residue degree 2, we first construct the unramified quadratic
extension of $\Q_5$.  By the procedure described in \S\ref{unram}, we
have $K^{u} = \Q_{5}[\alpha]/(\alpha^{2} - \alpha + 2)$.
Here $g=\gcd(e,p^f-1) = \gcd(3, 5^2-1) = 3$.  
Multiplication by $5$ on $\Z/3\Z$ has
two orbits, $\{1, 2\}$ and $\{0\}$, so there will be two extensions.
In the first case, $\{1, 5\}\subset \Z/24\Z$ is the prescribed lift,
so we take the norm of $x^3-5\alpha$ to get $x^6 - 5x^3 + 50$.
For the other orbit, the first orbit modulo 24 of length $f=2$
reducing to $\{0\}$ is $\{3, 15\}$.  Thus, we take the norm of
$x^3 -5 \alpha^3$ to get $x^6 + 25x^3 + 200$.

As an example where the root shift procedure is necessary, 
consider degree 12 extensions of $\Q_{5}$ with $e=6$ and $f=2$  
so that $g=\gcd(6,24)=6$.  The orbit $\{0\}\subset \Z/6\Z$
has only lifts of size 1 in $\Z/24\Z$.  However root shifting
with $k=1$ gives us the norm of $(x+\alpha)^6-5$, which is the
irreducible polynomial
$x^{12} +
6x^{11} + 27x^{10} + 80x^9 + 195x^8 + 366x^7 + 571x^6 + 702x^5 +
1005x^4 + 1140x^3 + 357x^2 - 138x + 44$.
%Luckily, the need for this
%process is uncommon among low degree extensions.
%

\subsection{Degree $p$ ramified extensions of $\Q_p$}
\label{sec:degree-p-ramified}
The six ramified quadratic extensions of $\Q_{2}$ are given 
by $x^{2}-D$ for $D = -4$, $12$, $\pm 8$, and $\pm 24$, with  
$\ord_{2}(D)$ being the discriminantal exponent $c$.  
Each of these six extensions has two automorphisms.  The rest 
of this subsection treats the case of $p$ odd, which is 
different as the generic degree $p$ extension of 
$\Q_{p}$ has just the identity automorphism.

%Since it is elementary to write down a defining equation
%for each of the six ramified quadratic extensions of
%$\Q_{2}$, 
%in this subsection we take $p$ odd. 
 Most of the
information we need can then be extracted from \cite{amano}.
Table~\ref{tab1} summarizes these results, giving exactly one
polynomial for each isomorphism class of degree $p$ extension of
$\Qp$.
\begin{table}[htb]
  \caption{\label{tab1}Degree $p$ ramified extensions of $\Q_p$, for 
  $p$ odd.}
  \begin{center}
  \begin{tabular}{|l|c|c|c|c|}
    \hline
%    \multicolumn{1}{|c|}{ \textbf{Family}} & \textbf{Parameters}
        \multicolumn{1}{|c|}{ \text{Family}} & \text{Parameters}
%    & \#
    & $c$ & $G$ & $I$ \\
    \hline
%
%    $x^p +apx^{p-1} + p$ & $0< a \le p-2 $
%    & $p-2$
%    & $2p-2$ & $C_p$
%    & $C_p.C_{d_2}$   \\
%    \hline
%
    $x^p +apx^\lambda + p$ &
    \begin{minipage}{1.40in}
      $\begin{array}{r@{}c@{}l}
    0< & a & \le p-1 \\
    1\le & \lambda & \le p-1 \\
    (\lambda,a) & \neq & (p-1,p-1) 
  \end{array}$
\end{minipage}
%    & $(p-1)(p-2)$
    & $p+\lambda-1$ & $C_p.C_{d_2}$
    & $C_p.C_{d_1}$   \\
%    & $1\le \lambda \le p-2$ & & & & \\
%
    \hline
    $x^p - p x^{p-1} + p(1+ap)$ & $0\le a \leq p-1$
%    & $p$
    & $2p-2$ & $C_p$  & $C_p$   \\
    \hline
%            $x^p+p(1+ap)$ & $0\le a \leq p-1$
%    & $p$
%    & $2p-1$    & $F_p$ & $F_p$   \\
%    \hline
                $x^p+p(1+ap)$ & $0\le a \leq p-1$
%    & $p$
    & $2p-1$    & $C_p.C_{p-1}$ & $C_{p}.C_{p-1}$   \\
    \hline
  \end{tabular}
  \end{center}
\end{table}
These come in three families, the main one with two parameters and 
the other two families each with one parameter.   
Table~\ref{tab1} gives our preferred defining
polynomials, restrictions on the parameter(s), the exponent $c$ of the
discriminant, and the Galois and inertia groups.  To define $d_1$ and
$d_2$, let $g=\gcd(p-1,c)$.  Then $d_1 = (p-1)/g$ and 
%$d_2 = \frac{p-1}{\gcd((p-1)/m,g)}$.
 $d_2 = (p-1)/(\gcd((p-1)/m,g))$
where
$m$ is the order of $a\lambda$ in $\F_p^*$.
%This applies as well to the
%third family by taking $\lambda = p-1$.

\subsection{Wild extensions of composite degree}
\label{sec:wild-extens-comp}
The complexity of the unramified, tamely ramified, and degree $p$ 
cases just treated suggests that analogous recipes for the remaining
cases would have to be quite complex indeed.  So instead, we treat the
five cases $(p,n) = (2,4)$, $(2,6)$, $(2,8)$, $(3,6)$, and
$(3,9)$ individually.  The problem
then becomes simply to find a defining polynomial for each degree $n$
extension of $\Qp$ for the given $(p,n)$. 

Pauli and Roblot give a general algorithm for
solving this problem.  One key ingredient is Panayi's  $p$-adic 
root finding algorithm \cite{panayi}, \cite[Section~8]{pauli-roblot}
which lets one determine whether two degree $n$
fields $\Q_{p}[x]/f_{1}(x)$ and $\Q_{p}[x]/f_{2}(x)$ are isomorphic 
and similarly lets one compute the set of automorphisms of a given
field $\Q_{p}[x]/f(x)$.  Another key ingredient is the
mass formula \cite[Theorem~6.1]{pauli-roblot} which lets one 
determine when all fields have been found.  

We used Pauli and Roblot's approach for generating polynomials as 
needed.  However, in many cases, we were able to generate 
the polynomials more
efficiently by specialized methods.  For example, to generate
2-adic fields of a given degree which have an index 2 subfield, we
were able to look up all possible candidates for subfields from our
database and construct the desired fields by taking square roots of
suitably chosen elements.  Similarly, when dealing with sextic
fields we worked with the twin algebra \cite{Ro-twin}. 
All except twelve 3-adic fields have 
reducible twin algebras. Thus, almost all sextic fields
could be generated by sextic twinning using lower degree
fields from the database.

%% a given which we use to
%% complete the lists.  The basic idea is to target fields of a given
%% degree and discriminant exponent, and to maintain a list of
%% polynomials defining non-isomorphic extensions.  Then one tests
%% randomly chosen polynomials from a collection sure to contain
%% representatives for each field in the target.  For example, for a
%% totally ramified extension, one searches for eisenstein polynomials
%% whose coefficients satisfy Ore-congruences to match the desired
%% discriminant exponent.  Each polynomial is tested against the current
%% list to see if it defines a new field.  One can tell when the list is
%% complete by means of mass formulae.

\section{Invariants associated to a given $p$-adic field}
\label{invariants}
Let $f(x) \in \Z[x]$ be a polynomial on one of our $p$-adic tables. 
In this section, we discuss the invariants  
the tables present for the corresponding field $K  = \Q_{p}[x]/f(x)$.  
Table~\ref{2adicquartics} serves as an example for much of our
discussion in this section.
\begin{table}[htb]
\caption{\label{2adicquartics}The first six lines of the 2-adic 
quartic table, corresponding to the fields with $c \leq 4$.}
{\renewcommand{\arraycolsep}{.05in}
\[
\begin{array}{|ccccc|c|cc|c|c|c|}
    \hline
 &  &  &  & &  &  &  & 
\mbox{Wild} &  & \mbox{Deg 2}  \\
c & e & f & \epsilon & d & \mbox{Polynomial} & G & I & 
\mbox{Slopes} & GMS & \mbox{Subs}    \\
\hline
0 & 1 & 4 & 1 & * & x^4+x+1 & C_{4} & <e> & & 0 & * \\
\hline
4 & 2 & 2 & -1 & 1 & x^{4} + 8 x^{2} + 4 & V_{4} & C_{2} & 2 & 1 & 
*,-1,-* \\
4 & 2 & 2& -1 & * & x^{4} - x^{2} + 5 & C_4 & C_{2} & 2 & 1 & * \\
4 & 2 & 2 & -i & -1 & x^{4} + 2 x^{2} + 4 x + 4 & D_{4} & V_{4} & 
2,2 & 3/2 & * \\
4 & 2 & 2 & -i & -* & x^{4} - 5 & D_{4} & V_{4} & 
2,2 & 3/2 & * \\
\hline
4 & 4 & 1 & 1 & * & x^{4} + 2 x + 2 & S_{4} & A_{4} & 4/3,4/3 & 7/6 &   
\\
\hline
\end{array}
\]
}
\end{table}

\subsection{Basic Data}
     The field discriminant of $K$ as an ideal is $(p^{c}) \subseteq \Z_{p}$.  
The largest unramified subfield of $K^{u}$ of $K$ has degree
the residual degree $f = [K^{u}:\Q_{p}]$.  The ramification
index is $e = n/f = [K:K^{u}]$.  The entry $d$ in the fifth column is the 
field discriminant considered as an element of 
$\Q_{p}^{\times}/\Q_{p}^{\times 2}$.  Here and elsewhere, 
$* \in \Q_{p}^{\times}/\Q_{p}^{\times 2}$ stands for a non-square unit.
With this notational convention, 
$\Q_{2}^{\times}/\Q_{2}^{\times 2} = \{1,*,-1,-*,2,2*,-2,-2*\}$
and otherwise  $\Q_{p}^{\times}/\Q_{p}^{\times 2} = \{1,*,p,p*\}$.
The computer program \gp{} \cite{gp} has commands to compute $c$, $e$, $f$, and 
$d$.

\subsection{Subfields and automorphisms}
For each field, we give its subfields hyperlinked to their respective
entries in the database.
Quadratic subfields are listed by the codes described in the previous
section.  An unramified subfield of degree $d$ is listed as simply
$U_d$.  All other subfields are listed by their chosen defining
polynomial.
To determine if one field is a subfield of
another, we make use of Panayi's $p$-adic root finding algorithm
mentioned in \S\ref{sec:wild-extens-comp}.  Similarly, we use this 
algorithm to find the automorphisms of $K$.    

%\cite{panayi, pauli-roblot}.  Given $f,g\in\Qp[x]$
%with $f$ irreducible, Panayi's algorithm efficiently determines if $g$
%has roots in $\Qp[x]/(f)$, and can be used to find approximations for
%those roots.  We have coded this algorithm into \gp.  In addition to 
%its used for subfield computations, is is used in
%\begin{itemize}
%\item the method of Pauli and Roblot described in
%  \S\ref{sec:wild-extens-comp}
%\item in taking an arbitrary polynomial which is irreducible over
%  $\Qp$ and finding its field in the database
%%\item in determining the automorphisms of a field, which is used in
%  the computation of its Galois group
%\end{itemize}

%% Moreover, for extensions of small degree, one can use Panayi's
%% algorithm with Algorithm~3.6.4 from \cite{cohen1} to factor
%% polynomials over finite extensions of $\Qp$ (as opposed to just
%% finding degree 1 factors).

%% We employ Panayi's algorithm to determine subfields of a given field.
%% Since we have complete lists of lower degree $p$-adic fields, each
%% field can be tested to see if its defining polynomial has a root in
%% the given field.  Subfields are listed by their defining polynomials
%% with two exceptions.  Quadratic subfields are coded in the same way
%% manner as discriminant root fields, and unramified subfields are
%% listed simply by $U_{m}$ where $m$ is the degree of the subfield.

\subsection{Root numbers}
%Root numbers are computed by,\ldots \textbf{hey, what is a reference
%for the root number computations?}
%
The root number $\epsilon$ is a complex fourth root of unity.  It is
of use in distinguishing fields, especially when $p=2$. 
However its principal use is in applications, for example to 
quadratic lifting, as explained in \cite[Section~2.1]{jr1}.

\subsection{Galois groups}
   Let $K^{g}$ be a Galois closure of $K$.  Our tables present
the isomorphism type of the local Galois group $G = \Gal(K^{g}/\Q_{p})$. 
Let $K^{g,u}$ be the maximum
unramified subfield of $K^g$.  
Our tables also present the inertia group 
$I = \Gal(K^{g}/K^{g,u})$, which is a normal subgroup
of $G$ such that $G/I$ is cyclic.  
Let $K^{g,t}$ be the maximal tame subfield of $K^g$.  Then 
the wild inertia group $I_{w} = \Gal(K^{g}/K^{g,t})$ 
is the normal subgroup of $I$ of $p$-power order
and prime-to-$p$ cyclic quotient.   
The quantity $t = |I/I_{w}|$, which appears in equation~(\ref{gms})
below, can then easily be read off.  

   The computation of Galois groups over $\Q_p$ is similar
to the more familiar computation of Galois groups over $\Q$.
For example, knowledge of the automorphism group 
$\mbox{Aut}(K/\Q_{p})$ and complete knowledge of subfields of 
$K$ restricts the possibilities.  Also factoring resolvents, now
over $\Qp$ of course, is the principal technique for distinguishing
between possible Galois groups.  However some techniques for computing
Galois groups over $\Q$ are not available when working over $\Qp$.  
For example, one cannot use cycle types of Frobenius elements 
for varying primes $p$.
In compensation, there are techniques which are particular to 
local fields.  Certainly extensions of $\Q_{p}$ are always solvable,
as indeed one has the wild-tame-unramified filtration
$I_{w} \le  I \le G$.  Also wild slopes, as discussed in the
next subsection, can be used to get lower bounds on the size of
$I_{w}$; in this sense they serve as substitute for Frobenius elements, 
which provide lower bounds in the global case.  

The discriminant class $d$ determines the parity of the
Galois group.  This much suffices for $n=3$, i.e.\ $G \cong A_{3}$ if
$d=1$ and $G \cong S_{3}$ if $d \neq 1$.  
Table~\ref{tab-low-galois} summarizes the computation in degrees $n=4$
and 5.  Here we use the number of automorphisms of 
the degree $n$ field $K$ to distinguish within  
$(C_{4},D_{4},S_{4})$, $(V_{4},A_{4})$,
and $(C_{5},D_{5})$.   
%Note that $A_5$ and $S_5$ are not candidates since they
%are not solvable.
Note that no resolvants are necessary here beyond using $d$ to
determine the parity of $G$, 
which is equivalent to considering the factorization of $x^{2}-d$.

\begin{table}[hbt]
  \caption{\label{tab-low-galois}Galois groups for $n=4$, $5$.}
  \begin{center}
  \begin{tabular}{|c|ccccc|ccc|}
    \hline
    $G$ & $C_4$ & $V_4$ & $D_4$ & $A_4$ & $S_4$ &
    $C_5$ & $D_5$ & $F_5$ \\ \hline
    Parity & $-$ &  $+$  & $-$ &  $+$  &$-$ & $+$ & $+$ & $-$\\ \hline
    $|\Aut(K)|$ & 4 & 4 & 2 & 1 & 1 & 5 & 1 & 1 \\ \hline
  \end{tabular}
  \end{center}
\end{table}

There are 16 transitive subgroups of $S_6$ up to conjugation, 12 of
which are solvable.  The Algorithm~6.3.10 of \cite{cohen1} starts by
computing a resolvant polynomial for the twin algebra of
$f$.  Most cases are determined by factoring this twin and considering
Galois groups of the irreducible factors, as shown by
Table~\ref{tab-deg6}. 
\cmmt{
\begin{table}[hbt]
  \caption{\label{tab-deg6}Galois groups for $n=6$.}
  \begin{tabular}{|c|*{6}{c}|}
    \hline
    $G$ & T1 & T2 & T3 & T4 & T5 & T6 \\ \hline
Twin & $C_3 C_2 C_1$ & $S_3 C_1^3$ & $S_3 C_2 C_1$ &
$A_4 C_1^2$ & $S_3 C_3$ & $A_4 C_2$ 
\\ \hline
\multicolumn{6}{c}{} \\ \hline
 $G$ & T7 & T8 & T9 & T10 & T11 & T13  \\ \hline
Twin & $S_4C_2$ & $S_4C_1^2$ & $S_3 S_3$ & Irred & $S_4 C_2$ &
Irred 
\\
    \hline
  \end{tabular}
\end{table}
}
\begin{table}[hbt]
  \caption{\label{tab-deg6}Galois groups for $n=6$.}
  \begin{center}
  \begin{tabular}{|c|*{6}{c}|}
    \hline
    $G$ & T1 & T2 & T3 & T4 & T5 & T6 \\ \hline
Twin & $C_3 C_2 C_1$ & $S_3 C_1^3$ & $S_3 C_2 C_1$ &
$A_4 C_1^2$ & $S_3 C_3$ & $A_4 C_2$ \\
\hline
Parity & - & - & - & + & - & - \\
\hline
$|\mbox{Aut}(K)|$ & 6 & 6 & 2 & 2 & 3 & 2 \\
 \hline
\multicolumn{6}{c}{} \\ \hline
 $G$ & T7 & T8 & T9 & T10 & T11 & T13  \\ \hline
Twin & $S_4^{+}C_2$ & $S_4C_1^2$ & $S_3 S_3$ & Irred & $S_4 C_2$ &
Irred \\
\hline
Parity & + & - & - & + & - & - \\
\hline
$|\mbox{Aut}(K)|$ & 2 & 2 & 1 & 1 & 2 & 1 \\
    \hline
  \end{tabular}
  \end{center}
\end{table}
Here, $C_1$ corresponds to linear factors in the twin
polynomial.
The two ambiguities remaining are distinguished by the parities of the
groups: T7 is even while T11 is odd, and
T10~$\cong
C_3^2.C_4$ is even while T13~$\cong C_3^2.D_4$ is odd.  Algorithm~6.3.10
mentioned above sometimes requires the computation of a second
resolvant.  Here that is unnecessary because of the absence of
non-solvable groups.  Also here the information on $|\Aut(K)|$ is
not needed to distinguish Galois groups; we present it because of
the important role it plays in the mass formula.  

In $S_7$, there are 7 transitive subgroups of which 4 are
solvable, namely $C_{7}.C_{d}$ for $d=1$, $2$, $3$, and $6$.
One could apply Algorithm~6.3.11 of \cite{cohen1} directly.  It
uses the discriminant and the factorization of a degree 35 resolvant.
However, this high degree resolvant can be avoided since all cases are
covered by Section~\ref{construction}.
%as follows.  
If $K$ is unramified, then $G = C_{7}$.  If $G$ is tamely ramified, 
then $G = C_{7}.C_{d}$ with $d$ the order of the subgroup 
of $\F_{7}^{\times}$ generated by $p$.  If $G$ is wildly ramified, then
Table~\ref{tab1} covers the situation.
% New paragraph?
%
In the much more complicated 
cases $n=8$ and $n=9$,  we similarly use specifically $p$-adic 
phenomena to avoid large degree resolvents.  

%Ironically, in degrees 6 and 7, our approach starts with a polynomial
%viewed over $\Qp$, yet the resolvant computation makes use of its
%complex roots.  Our polynomials have integer
%coefficients, so when we compute resolvant polynomials, we employ the
%standard approach of using the polynomial's complex roots computed to
%high precision, constructing the resolvant, and rounding off the
%coefficients of the resulting polynomials.  The resolvant is, of
%course, the same as if it had been computed using roots in
%$\overline{\Q}_p$, so we can proceed and factor it over $\Qp$.

%Computing Galois groups of $p$-adic fields in degrees 8 and 9 are
%sufficiently involved that we will discuss them in \cite{jr-octic} and
%\cite{jr-nonic}.  The same comment applies to the computations in the
%remainder of this section, where we will restrict our attention to
%degrees of at most $7$.

\subsection{Slope filtration and root discriminants}
For each field we compute its slope content as in \cite{jr1},
whose definition we briefly recall here.  Let $G^\nu$ be the ramification
groups with Artin upper numbering so that $G^{0}=G$ and 
$G^1 = I$.  Non-trivial quotients $Q^\nu = G^{\nu}/G^{>\nu}$ contribute
slope $\nu$.  In terms of the wild-tame-unramified filtration
mentioned above, $Q^0 = G/I$ corresponds to the unramified part, $Q^1
= I/I_w$ corresponds to the tame part, and slopes greater than 1
contribute to the wild part.  
We have discussed finding $|Q^{0}|$ and $|Q^{1}|$ in the previous
subsection, so our focus now is finding the wild slopes $\nu$.
Our tables list each wild slope $\nu$, repeated according
to its multiplicity $m = \ord_{p}(|Q^\nu|)$.
So the total number of wild slopes listed is $\ord_{p}(|I|) = 
\log_{p}(|I_{w}|)$.   
%For convenience, when speaking of multiplicity we also 
%consider $0$ to have multiplicity $\ord_{p}(|Q^{0}|)$.  The the total number of $p$-slopes is
%$\ord_{p}(|G|)$.  

%Tame and unramified parts always have slopes 1 and 0 respectively -
%what is of interest there are their degrees.  In particular, we let
%$t=|Q^1|$.  Wild slopes, on the other hand, always correspond to
%subquotients of order $p$; the slope (counted with multiplicity) is of
%interest rather than the corresponding degrees.

%Since one may have $p$ dividing $|Q^0|$, we refer to $p$ having
%slopes of 0, the multiplicity again referring to the exponent of $p$ in
%$|Q^0|$.  Thus, the total number of $p$-slopes associated to an
%extension is the exponent of $p$ in $|G|$.

To compute slopes for a field, we use the following proposition.
\begin{Proposition} \label{slope-tower}
Let $K$ be a $p$-adic field.  Say a subfield is distinguished iff all 
other subfields of the same degree have strictly larger 
discriminant exponent.  Then the distinguished subfields form a chain
\[
\Q_{p} = K_{0} \subset K_{1} \subset \cdots K_{n-1} \subset K_{n} = K.
\]
Say $K_{i}$ has degree $n_{i}$ and discriminant $p^{c_{i}}$.  Then the 
\begin{equation}
\label{slopeformula}
s_{i} = \frac{c_{i} - c_{i-1}}{n_{i} - n_{i-1}} 
\end{equation}
are weakly increasing.  If $s_{j+1}=s_{j+2}= \cdots = s_{i}$, then
$|Q^{s_{i}}| \geq n_{i}/n_{j}$.  
\end{Proposition}
\noindent This proposition follows from the basic facts of ramification theory
as follows.  An irreducible representation
\[\rho\colon G \to \Aut(V) \]
of $G$ has a slope $s(\rho)$, namely the smallest $c \in [0,\infty)$ with
$G^{>c}$ in the kernel of $\rho$.  The Artin exponent $c(\rho)$ 
of $\rho$ is the slope $s(\rho)$ times the degree $\dim(V)$.
Artin exponents of arbitrary representations behave additively. 
If $\rho$ is induced from
a permutation representation of $G$ on  a finite set $X$, 
%$r\colon G \to S_X$,
then the discriminant of the $p$-adic algebra corresponding 
to $X$ is $p^{c(\rho)}$.  Take $X_{i}$ corresponding to $K_{i}$
so that the inclusions $K_{i-1} \subset K_{i}$ give 
surjections $X_{i} \rightarrow X_{i-1}$.  The action
of $G$ on $V = \Q[X_{i}]/\Q[X_{i-1}]$ must have a single slope
as otherwise there would be a distinguished subfield between
$K_{i-1}$ and $K_{i}$.  As the Artin exponent of $V$ is
$c_{i} - c_{i-1}$ and the dimension of $V$ is $n_{i} - n_{i-1}$, one
gets (\ref{slopeformula}).

To analyze a field $K$, we begin by applying 
Proposition~\ref{slope-tower} directly to $K$.  We refer to the 
wild slopes we see here as visible slopes.  In many low degree cases,
these suffice.  For example, for a ramified degree $p$ extension of $\Qp$, $p$
exactly divides the order of the inertia group.  Thus, there is exactly
one wild slope and it is visible from the degree $p$ extension, being
just $c/(p-1)$.  When we have not found enough slopes to account for 
all of $G$, we apply Proposition~\ref{slope-tower} to various
resolvent fields, typically the same resolvents used to compute
$G$.  

For example, the first case of wild ramification with composite degree
is $2$-adic quartics.  Here up to two wild slopes will be visible from
the quartic field.
When the Galois group is $C_4$, $V_{4}$, $A_4$, or $S_{4}$, this 
suffices, the visible slopes being the only slopes.   
The remaining case for $2$-adic quartics is when the Galois group is
$D_4$, where there are several viable approaches.  One would be to
compute the octic Galois field by taking the compositum of the defining
quartic and $x^2-d$.  At this octic level, all slopes are visible.  
% For example, there are twelve $2$-adic quartics with $G = D_{4}$.  
% All have two visible wild slopes, namely $3$ and $4$.  For four of 
% them, the slopes are $0$, $3$, and $4$ while for the other eight
% the slopes are $2$, $3$, and $4$.  

When computing slopes of sextic fields, we make use of the sextic twin
algebra, which we compute as part of the Galois group
computation.  In most cases, this is the product of smaller degree
fields.  On the one hand, we will have already have computed the slopes
for these smaller fields.  On the other hand, when considering the
composita of the lower degree fields, slopes may combine in
non-trivial ways, as discussed below in
\S\ref{sec:grd-computations}.  

The twin algebra of a solvable sextic field is also a field
when the Galois group of the normal closure is $C_3^2.C_4$ or
$C_3^2.D_4$.   These only appear for $p=3$, where the first appears
four times and the second eight times.  In both cases, the chain of 
distinguished fields takes the form $\Q_{3} \subset K_{2} \subset 
K_{6}$, with the latter two fields now indexed by degree. 
Proposition~\ref{slope-tower} then says that $(c_6-c_2)/4$ is 
a wild slope.  Since neither group contains a
normal subgroup of order 3, this slope must be repeated with
multiplicity 2.  This is an instance where group theory makes 
resolvent constructions unnecessary.  The cases of degrees eight and
nine involved group-theoretical arguments of this type, as well
as actual resolvent constructions.  

\subsection{Galois mean slope}  The root discriminant of $K$ is
$p^{c/n}$.  Also interesting is the corresponding quantity for the 
Galois closure $K^{g}$, which we write as $p^{\beta}$ with 
$\beta$ what we call the Galois mean slope.  This Galois mean slope can be computed 
as a weighted sum of the slopes, with larger slopes counted more 
heavily.   More precisely, let $s_{1},\dots,s_{m}$ be the slopes in
decreasing order,  so that $|I_{w}|=p^{m}$.  Then 
\begin{equation}
\label{gms}
\beta = \left(  \sum_{i=1}^{m} \frac{p-1}{p^{i}} s_{i} \right) + 
\frac{1}{p^{m}} \frac{t-1}{t},
\end{equation}
with $t = |G/I|$.  
One has $\beta \geq c/n$ with equality iff $K^{g} \cong K \otimes U$
with $U$ an unramified extension of $\Q_{p}$, i.e.\ iff there are
no hidden slopes.  The Galois mean slopes $\beta$ play an important 
role in \S\ref{sec:grd-computations}.

\section{A sample use of the database}
\label{sampleuse}
Our database has two interactive features, the {\em $p$-adic identifier} and 
the {\em GRD calculator}.  Here we illustrate how they can be used in 
the study of number fields by means of a family of examples.  
\subsection{A family of number fields}
Let
\begin{eqnarray*}
\lefteqn{f_{t}(x) =} \\   
&& (x^{6} + 12 x^{5} + 54 x^{4} + 176 x^{3} + 444 x^{2} + 624 
x + 552)  \cdot \\ &&(x^{4} + 16 x^{3} - 36 x^{2} + 128 x - 28)(x^{3} - 12 
x^{2} -  
6 x - 64) \\
& & + t (3 x^{4} - 4 x^{3} + 12 x^{2} - 24 x -68)^{2} (x^{4} + 16 
x^{3} + 72 x^{2} + 128 x + 188).
\end{eqnarray*}
This one parameter family of polynomials was first introduced in \cite{malle-jsc}.
Its polynomial 
discriminant is $2^{160} 3^{114} (t^{2}+8)^{10}$.  We are 
interested in the fields $K_{t} = \Q[x]/f_{t}(x)$ for $t \in \Q$.  
Let $K^{g}_{t}$ be the splitting field of $f_{t}(x)$ in $\C$.  
For generic $t$, the Galois group $\Gal(K_{t}^{g}/\Q)$ is isomorphic 
to the projective linear group $PSL_{3}(\F_{3})$, of order 
$2^{4}3^{3}13 = 5616$.  For these $t$, $K_{-t} = \Q[x]/f_{-t}(x)$ is 
the projective twin of $K_{t}$, meaning that if $K_{t}$ corresponds
to an action of $\Gal(\Qbar/\Q)$ on
the projective space ${\bbP}^{2}(\F_{3})$, then 
$K_{-t}$ corresponds to the induced action on the dual projective 
space $\check{\bbP}^{2}(\F_{3})$.  So $K_{t}$ and $K_{-t}$ are 
non-isomorphic fields with the same Dedekind zeta function, the 
same discriminant, and the same Galois closure $K_{t}^{g} = 
K_{-t}^{g}$.  If $p$ is such that the local decomposition group 
contains an element of order $13$ then both 
$K_{t} \otimes \Q_{p}$ and $K_{-t} \otimes \Q_{p}$ are degree 13 
fields, thus treated by the dynamic part of our database.  Otherwise,
by consideration of maximal subgroups of $PSL_{3}(\F_3)$, at least one
of $K_{t}$ and $K_{-t}$ splits $p$-adically into a nonic times a 
quartic, and perhaps further.       

\subsection{The $p$-adic identifier}
    The $p$-adic identifier lets one input a polynomial 
$f(x) \in  \Z[x]$ and a prime $p$.  
It uses Panayi's root finding algorithm and 
returns the entries from the database corresponding to the
factor fields of $\Q_{p}[x]/f(x)$.  In the context of the family 
$f_{t}(x)$, it lets one study the $p$-adic completion
of these fields, as a function of the parameter $t$, now
considered $p$-adically.  The interesting cases are $p=2$ and
$p=3$, as otherwise ramification is tame.  For $p=2$, the
behavior is locally constant for $t \in \Q_{2}$.  Experiment
quickly suggests that the factorization pattern is $4+4+2+2+1$
if $\ord_{2}(t) \geq 4$ and otherwise $8+4+1$, with many possible
fields arising for the octic.  For $p=3$, the behavior is locally
constant for $t \in \Q_{3} - \{\sqrt{-8},-\sqrt{-8}\}$ with
$\sqrt{-8} = 1 +  3^{2} + 3^{3} + 2 \cdot 3^{4} + 3^{5} + 3^{6} + 2 
\cdot 3^{8} + \cdots \approx 14293$.  Here again the fields involved
vary substantially with $t$.  However experiment
suggests that factor fields are always of degree $\leq 9$, unless
$\ord_{3}(t - \sqrt{-8}) \geq 3$, in which case they may be
of degree $12$.

\begin{table}[htb]
\caption{\label{identifiertable}Information 
output by the $p$-adic identifier for the polynomial $g_{-7/2}(x)$.  
The field factors $2$-adically as $\Q_{2} \times K_{4} \times K_{8}$ 
with information on $K_{4}$ and $K_{8}$ being given in the first block.
The field factors $3$-adically as $K_{4} \times K_{9}$
with information on $K_{4}$ and $K_{9}$ given in the second block.}
\begin{center}
{\renewcommand{\arraycolsep}{.04in}
\[
\begin{array}{|ccccc|c|cc|c|c|c|}
    \hline
 &  &  &  &  & &  &  & 
\mbox{Wild} & &  \\
c & e & f & \epsilon & d & \mbox{Polynomial} & G & I & 
\mbox{Slopes} & GMS  & \mbox{Subfields}  \\
\hline
\multicolumn{11}{c}{\;} \\
\hline
6 & 2 & 2 & i & -1 & x^{4} + 2 x^{2} - 4 & D_{4} & V_{4} & 2,3 & 
2 & * \\
22 & 4 & 2 & -i & -1 & x^{8} + 10 x^{4} + 28 & D_{8} & D_{4} & 2,3,4 & 
3 & *,x^{4}-2 x^{2} + 4 \\
\hline
\multicolumn{11}{c}{\;} \\
\hline
3 & 4 & 1 & -i & 3 & x^{4}+3 & D_{4} & C_{4} &   & 3/4 & 3*  \\
19 & 9 & 1 & i & 3 & x^{9}+9 x^{2}+3& T19 & T15 & 19/8,19/8 & 53/24 &  \\
\hline
\end{array}
\]
}
\end{center}
\end{table}

   As a specific example, consider $t = \pm 7/2$, which we choose 
because the corresponding fields have the smallest discriminant we 
are aware of for $PSL_{2}(\F_3)$ fields, namely $2^{28} 3^{22}$.  
For comparison, the database \cite{kluner-malle} has $2^{28} 3^{24}$ as its 
smallest discriminant, and the two fields given there
are $K_{\pm 1}$.  Applying \gp's {\em polredabs} to $f_{\pm 7/2}(x)$ to get 
smaller coefficients, we have
\begin{eqnarray*}
g_{-7/2}(x) & = &
x^{13} - 2 x^{12} - 8 x^{10} + 55 x^{9} - 90 x^{8} -108 x^{7} + \\
&&
684 x^{6} -1341 x^{5} + 1526 x^{4} - 1090 x^{3} + 468 x^{2} - 100 x + 8 \\
g_{7/2}(x) & = &
x^{13} - x^{12} - 3 x^{11} - 7 x^{10} + 37 x^{9} - 9 x^{8} - 168 x^{7} + \\
&& 24 x^{6} + 396 x^{5} + 20 x^{4} - 128 x^{3} + 192 x^{2} - 176 x - 16.
\end{eqnarray*}
Entering $f_{-7/2}(x)$ or equally well $g_{-7/2}(x)$, one gets the
$2$-adic information listed in Table~\ref{identifiertable}.  One also 
gets 
that the polynomial factors 3-adically as 12+1, not a surprise, as 
\[
\ord_{3}(-7/2-\sqrt{-8}) = \ord_{3}(-7/2 - 14293) = \ord_{3}(-3^{4} 
353/2) = 4. 
\]
Entering instead $f_{7/2}(x)$ or $g_{7/2}(x)$, one 
gets the same identification at 2 and now a 3-adic factorization of
9+4, with information as on Table~\ref{identifiertable}.

\subsection{The GRD calculator}
\label{sec:grd-computations} 
A single numerical measure of ramification in a polynomial $f(x) \in 
\Z[x]$ is the 
root discriminant of its splitting field in $\C$, 
what we call its Galois root discriminant.  The GRD calculator accepts a 
polynomial 
$f(x) \in \Z[x]$ as input.  When all factors of all completions 
of $\Q[x]/f(x)$ are in the database, it returns lower and upper bounds 
on the Galois mean slope $\beta_{p}$
of each ramifying prime $p$, and hence bounds
on the Galois root discriminant $\prod p^{\beta_{p}}$.  In favorable
cases, certainly when the $p$-adic algebra has only one wild factor, 
the lower and upper bounds on $\beta_{p}$ agree.  In the remaining
cases, it is typically easy to start with the bounds and continue
by hand to exactly determine $\beta_{p}$, as our example will 
illustrate.  

The reason the GRD calculator returns only bounds is as follows.
We are assuming that all the $p$-adic factors $K_{p,i}$ of
$K \otimes \Q_{p}$ are in our database.  Thus we know the
corresponding tame parts of inertia $t_{i}$ and the multiplicities
$m_{i}(s)$ of any given wild slope $s>1$.
We need the corresponding information $t$ and $m(s)$ 
for the algebra $\prod_{i} K_{p,i}$.  The tame index for the algebra is simply 
$t = \lcm(t_{i})$, and this formula is incorporated into the GRD 
calculator.    The biggest possible Galois mean slope would arise if 
$m(s) = \sum_{i} m_{i}(s)$ for all $s$. 
The smallest possible Galois mean slope 
would arise if $m(s) = \max_{i} m_{i}(s)$ for all $s$.  
When these bounds disagree, i.e. when there is some overlap
between the factor wild slope lists, there may be quite a number of 
possibilities in between, including slopes in the algebra which 
are not in any of the factors.  For example, suppose $K_{p,1}$ and
$K_{p,2}$ are distinct degree $p$ fields each with $c=2p-1$,
thus each with unique wild slope $(2p-1)/(p-1)$.  
Then the slopes of $K_{p,1} \times K_{p,2}$ are $(2p-1)/(p-1)$
and $p/(p-1)$.  

The point of symmetry $t=0$ of our family of polynomials is forced to 
be special with respect to Galois theory as $K_{0}$ is its own 
projective twin.  It factors over $\Q$ as 6+4+3 and its 
Galois group is $S_{4}$.  The small size of the global Galois group
limits ramification at $2$ and $3$, and the GRD is the relatively 
small number $2^{11/4} 3^{7/8} \approx 24.2367$.  All this is recorded
as the top entry of Table~\ref{GRDtable}.

\begin{table}[htb]
\caption{\label{GRDtable}The seven smallest GRDs found in the
family $f_{t}(x)$.  The first five correspond to globally reducible
degree $13$ polynomials, as indicated by the second block of columns. 
The last two correspond to irreducible degree $13$ polynomials.}
\[
\begin{array}{|c|cc|rcr|r|c|}
\hline
t & K_{-t} & K_{t} & \multicolumn{3}{c|}{GRD} & |G| & \ord_{3}(t^{2}+8)\\
\hline
0 & \multicolumn{2}{c|}{6,4,3} & 2^{11/4} 3^{7/6} & \approx &  
24.2367 & 24 & 0 \\
23/5 & 12,1 & 9,4 & 2^{11/4} 3^{7/8} 5^{1/2}  & \approx & 39.3368 & 
144 & 6\\ 
10 & 9,4 & 12,1 & 2^{29/8} 3^{7/6} & \approx & 44.4504  & 432 & 3 \\
6808/2209 & 12,1 & 9,4 & 2^{3} 3^{7/8} 11^{1/2} & \approx & 69.3853 & 
144 & 6 \\
10516/725 & 12,1 & 9,4 & 2^{29/8} 29^{1/2} & \approx & 66.4405 & 432 & 
15 \hspace{5pt} \\  
7/2 & 13 & 13 & 2^{3} 3^{53/24} & \approx & 90.5175 & 5616 &  4 \\
1 & 13 & 13 & 2^{11/4} 3^{133/54} & \approx & 100.6849 & 5616 & 2 \\
\hline
\end{array}
\]
\end{table}

We have made a substantial search over carefully chosen rational $t$ and 
Table~\ref{GRDtable} presents the $t$ we have found with GRD less than 
110.  In each case, the lower bound returned by the GRD calculator 
is exact.  
%In each case, exactness can easily be verified by 
%hand. 
For example, consider the 2-adic Galois mean slope in the
case $t = \pm 7/2$, so that Table~\ref{identifiertable} applies.
A local reason that one need not consider the quartic and its contribution
of $2,3$ is that the quartic is a subfield of the octic.
We plan to incorporate this sort of argument into 
the GRD calculator in the future, although many local arguments of
this sort are more subtle.  A global reason that 
one need not consider the quartic's contribution is that the octic's
Galois group already has $16$ elements and $16$ exactly divides
$5616 = 2^{4} 3^{3} 13$.  Variants on this sort 
of global argument work in many situations.

Table~\ref{GRDtable} shows that for the four pairs $\pm t$ yielding
the smallest GRD, the Galois group drops down from its generic
size of 5616.  The column $\ord_{3}(t^{2}+8)$ gives the
maximum of the numbers $\ord_{3}(t+\sqrt{-8})$ and 
$\ord_{3}(t-\sqrt{-8})$, the other being zero.  Experiment suggests 
that ramification can be tame at $3$ only if this number is $\geq 6$ 
and $3$ can be unramified only if this number is $\geq 15$, while
$2$ is always wildly ramified with Galois mean slope at 
least $11/4$; these observations guided our search for $t$.
Many similar computations are presented in 
\cite{jrlowgrd}, where the ``low GRD Galois drop'' phenomenon visible in 
Table~\ref{GRDtable} is evident in similar strength.  Detailed global 
analyses such as these would not be possible without quite complete
control over local phenomena.

%%%%%%%%%%%%%%%%%%%%%%%%%%% Bibliography %%%%%

\begin{comment}

\end{comment}

\bibliography{jr}

\end{document}